\documentclass[12pt]{amsart}

\usepackage{amsmath}
\usepackage{amssymb}
\usepackage{array}

\renewcommand{\atop}[2]{%
\genfrac{}{}{0pt}{}{#1}{#2}}

\newtheorem{theorem}{Theorem}

\theoremstyle{definition}

\renewcommand{\atop}[2]{%
\genfrac{}{}{0pt}{}{#1}{#2}}

\begin{document}

\title{A $q$-analog of the Bailey-Borwein-Bradley identity}

\author{Kh.~Hessami Pilehrood}

\address{Mathematics Department, Faculty of Basic  Sciences,  Shahrekord
University,  Shahrekord, P.O. Box 115, Iran.}
\email{hessamik@ipm.ir, hessamit@ipm.ir, hessamit@gmail.com}

\author{T.~Hessami Pilehrood}


\date{}


\begin{abstract}
We establish a $q$-analogue of the Bailey-Borwein-Bradley identity
generating accelerated series for even zeta values and prove
$q$-analogues of Markov's and Amdeberhan's series for $\zeta(3)$
using the $q$-Markov-WZ method.
\end{abstract}

\maketitle

\section{Introduction}
\label{intro}

The Riemann zeta function for ${\rm Re}(s)>1$ is defined by the
series
$$
\zeta(s)=\sum_{n=1}^{\infty}\frac{1}{n^s}={}\sb {s+1}F\sb
s\left(\left.\atop{1,\ldots,1}{2,\ldots,2}\right|1\right),
$$
where
$$
{}\sb rF\sb t\left(\left.\atop{a_1,\ldots,
a_r}{b_1,\ldots,b_t}\right|z\right)=
\sum_{n=0}^{\infty}\frac{(a_1)_n\cdots (a_r)_n}{(b_1)_n\cdots
(b_t)_n}\frac{z^n}{n!}
$$
is the generalized hypergeometric function and $(a)_n$ is the
shifted factorial defined by $(a)_n=a(a+1)\cdots (a+n-1),$ $n\ge
1,$ and $(a)_0=1.$

 Ap\'ery's irrationality proof of
$\zeta(3)$ \cite{po} operates with the faster convergent series
\begin{equation}
\zeta(3)=\frac{5}{2}\sum_{k=1}^{\infty}\frac{(-1)^{k-1}}{k^3\binom{2k}{k}}
\label{z3}
\end{equation}
 first obtained by A.~A.~Markov in 1890 \cite{ma}.
Since then, many Ap\'ery-like formulas for other values
$\zeta(n),$ $n\ge 2,$ have been proved with the help of generating
function identities. In 1979 Koecher \cite{ko} (and independently
 Leshchiner \cite{le}) proved that
\begin{equation}
\sum_{s=0}^{\infty}\zeta(2s+3)x^{2s}=\sum_{k=1}^{\infty}\frac{1}{k(k^2-x^2)}=
\frac{1}{2}\sum_{k=1}^{\infty}\frac{(-1)^{k-1}}{k^3\binom{2k}{k}}\,\,\frac{5k^2-x^2}%
{k^2-x^2}\,\prod_{m=1}^{k-1}\left(1-\frac{x^2}{m^2}\right),
\label{eq01}
\end{equation}
for any $x\in {\mathbb C},$ with $|x|<1.$ Expanding the right-hand
side of (\ref{eq01}) in powers of $x^2$ and comparing coefficients
of $x^{2s}$ on both sides of (\ref{eq01}) gives Ap\'ery-like
series for $\zeta(2s+3).$  The similar identity generating fast
convergent series for all $\zeta(4s+3),$ $s\ge 0,$ which for $s>1$
are different from Koecher's result (\ref{eq01}) was
experimentally discovered in \cite{bb} and proved by G.~Almkvist
and A.~Granville in \cite{algr}
\begin{equation}
\sum_{s=0}^{\infty}\zeta(4s+3)x^{4s}=\sum_{k=1}^{\infty}\frac{k}%
{k^4-x^4}=\frac{5}{2}\sum_{k=1}^{\infty}\frac{(-1)^{k-1}}%
{\binom{2k}{k}}\frac{k}{k^4-x^4}\prod_{m=1}^{k-1}\left(\frac{m^4+4x^4}%
{m^4-x^4}\right). \label{eq02}
\end{equation}
There exists a bivariate unifying formula for identities
(\ref{eq01}) and (\ref{eq02})
\begin{equation}
\sum_{k=1}^{\infty}\frac{k}{k^4-x^2k^2-y^4}=\frac{1}{2}\sum_{k=1}^{\infty}
\frac{(-1)^{k-1}}{k\binom{2k}{k}}\frac{5k^2-x^2}{k^4-x^2k^2-y^4}\prod_{m=1}^{k-1}
\left(\frac{(m^2-x^2)^2+4y^4}{m^4-x^2m^2-y^4}\right), \label{eq03}
\end{equation}
which was first conjectured by H.~Cohen and then proved by
D.~Bradley \cite{br} and, independently, by T.~Rivoal \cite{ri}.
This identity implies (\ref{eq01}) if $y=0,$ and gives
(\ref{eq02}) if $x=0.$

For even zeta values, D.~Bailey, J.~Borwein and D.~Bradley
\cite{bbb} showed that
\begin{equation}
\sum_{k=0}^{\infty}\zeta(2k+2)a^{2k}=\sum_{n=1}^{\infty}\frac{1}{n^2-a^2}=
3\sum_{k=1}^{\infty}\frac{1}{\binom{2k}{k}(k^2-a^2)}\,\prod_{m=1}^{k-1}
\left(\frac{m^2-4a^2}{m^2-a^2}\right), \label{eq04}
\end{equation}
 for any $a\in
{\mathbb C}, |a|<1.$ In particular, if $a=0$ we get an
Ap\'ery-like formula for $\zeta(2)$
$$
\zeta(2)=\sum_{n=1}^{\infty}\frac{1}{n^2}=
3\sum_{k=1}^{\infty}\frac{1}{k^2\binom{2k}{k}}.
$$
In \cite{he1, he2} it  was shown that identities
(\ref{eq01})--(\ref{eq04}) and many other generating function
identities producing fast convergent series for zeta values can be
proved with the help of the Markov-Wilf-Zeilberger theory.

We say that a function $a(n)$ is {\it P-recursive} if it satisfies
a linear recurrence relation with polynomial coefficients. A
function $H(n,k)$, in the integer variables $n$ and $k,$ is called
a {\it hypergeometric term}  if the quotients
$$
\frac{H(n+1,k)}{H(n,k)} \qquad\mbox{and} \qquad
\frac{H(n,k+1)}{H(n,k)}
$$
are both rational functions of $n$ and $k.$
 If for a given
hypergeometric term $H(n,k),$  there exists a polynomial $P(n,k)$
in $k$ of the form
$$
P(n,k)=a_0(n)+a_1(n)k+\cdots+a_L(n)k^L,
$$
for some non-negative integer $L,$ and P-recursive functions
$a_0(n), \ldots, a_L(n)$ such that $F(n,k):=H(n,k)P(n,k)$
satisfies
\begin{equation}
F(n+1,k)-F(n,k)=G(n,k+1)-G(n,k) \label{WZ}
\end{equation}
with some function $G,$ then a pair $(F,G)$ is called a {\it
Markov-WZ pair} associated with the kernel $H(n,k)$ (MWZ-pair for
short). 

Considering the Markov-WZ pair \cite[\S 2]{he1}
$$
F(n,k)=(-1)^n (1+x)_n(1-x)_n H(n,k),
$$
$$
G(n,k)=H(n,k)\frac{(-1)^n(1+x)_n(1-x)_n}%
{(2n+k+2)(2n+2)}(5(n+1)^2-x^2+k^2+4k(n+1))
$$
associated with the kernel
$$
H(n,k)=\frac{k!}{(2n+k+1)!((n+k+1)^2-x^2)}
$$
proves Koecher's identity.

The following MWZ-pair \cite[\S 4]{he1}:
$$
F(n,k)=H(n,k)\frac{n!^2}{(2n)!}(1+2a)_n(1-2a)_n,
$$
$$
G(n,k)=H(n,k)\frac{n!(n+1)!}{(2n+2)!}(1+2a)_n(1-2a)_n(3n+3+2k)
$$
associated with the kernel
$$
H(n,k)=\frac{(1+a)_k(1-a)_k}{(1+a)_{n+k+1}(1-a)_{n+k+1}}
$$
gives a proof of the Bailey-Borwein-Bradley identity.

If we consider the Markov-WZ pair  (see \cite{he2})
$$
F(n,k)=\frac{H(n,k)}{2}\frac{(-1)^n}{\binom{2n}{n}} (1\pm a\pm
b)_n(n+2+2k),
$$
$$
G(n,k)=\frac{H(n,k)}{2}\frac{(-1)^n}{\binom{2n}{n}}\frac{(1\pm a
\pm b)_n}{2n+1}(5(n+1)^2-a^2-b^2+3(n+1)k+k^2)
$$
associated with the kernel
$$
H(n,k)=\frac{(1\pm a)_k(1\pm b)_k}{(1\pm a)_{n+k+1}(1\pm
b)_{n+k+1}}
$$
we derive the following identity:
$$
\sum_{k=1}^{\infty}\frac{k}{(k^2-a^2)(k^2-b^2)}=\frac{1}{2}
\sum_{n=1}^{\infty}\frac{(-1)^{n-1}(5n^2-a^2-b^2)(1\pm a\pm
b)_{n-1}}{n\binom{2n}{n}(1\pm a)_n(1\pm b)_n}.
$$
Here and below $(u\pm v\pm w)$ means that the product contains the
factors $u+v+w, u+v-w, u-v+w, u-v-w.$
 If we now put
$$
 a^2=\frac{x^2+\sqrt{x^4+4y^4}}{2},
\qquad b^2=\frac{x^2-\sqrt{x^4+4y^4}}{2},
$$
 we get the bivariate
identity (\ref{eq03}) conjectured by H.~Cohen.

The aim of the present paper is to establish a $q$-analogue of the
Bailey-Borwein-Bradley identity and to prove a family of
$q$-analogues of Markov's formula (\ref{z3}) for $\zeta(3)$ using
the $q$-Markov-WZ method. The formulas found imply that the
$q$-extensions of generating function identities for odd zeta
values (\ref{eq01})--(\ref{eq03}), if they exist, are likely to be
rather complicated.

\section{$q$-analogues of zeta values} \label{S2}

Throughout the paper unless otherwise stated we assume that
$0<q<1.$ The $q$-analogue of a positive integer $n$ is
$$
[n]_q=\frac{1-q^n}{1-q}=\sum_{k=0}^{n-1}q^k,
$$
the $q$-shifted factorial is defined by the product
\begin{equation*}
(a;q)_n =\begin{cases}
 1,     & \quad n=0; \\
(1-a)(1-aq)\cdots (1-aq^{n-1}), & \quad n\in {\mathbb N},
\end{cases}
\end{equation*}
and the $q$-binomial coefficient is given by the ratio
$$
\genfrac{[}{]}{0pt}{}{n}{k}_q=\frac{(q;q)_n}{(q;q)_k(q;q)_{n-k}}
\qquad n,k\in {\mathbb Z}, \quad n\ge k\ge 0.
$$
We define a $q$-analogue of the Riemann zeta function following to
the work of Kaneko, Kurokawa and Wakayama \cite{kkw}
$$
\zeta[s]=\sum_{n=1}^{\infty}\frac{q^{n(s-1)}}{[n]_q^s}.
$$
This definition is natural in view of the equality
$$
\zeta[s]=\sum_{n=0}^{\infty}\frac{(q;q)_n^s}{(q^2;q)_n^s}q^{(n+1)(s-1)}
=q^{s-1}{}\sb {s+1}\phi\sb s\left(\atop{q,\ldots,
q}{q^2,\ldots,q^2}; q, q^{s-1}\right),
$$
where
$$
{}\sb r\phi\sb t\left(\atop{a_1,\ldots, a_r}{b_1,\ldots,b_t}; q,
z\right)= \sum_{n=0}^{\infty}\frac{(a_1;q)_n\cdots
(a_r;q)_n}{(b_1;q)_n\cdots
(b_t;q)_n(q;q)_n}z^n((-1)^nq^{n(n-1)/2})^{1+t-r}
$$
is the basic hypergeometric series. The ordinary hypergeometry is
a limiting case of the basic one:
$$
\lim_{q\to 1}\frac{(q^a;q)_n}{(q^b;q)_n}=\frac{(a)_n}{(b)_n}.
$$
Therefore, we get for ${\rm Re}\, s>1$
$$
\lim_{q\to 1}\zeta[s]=\zeta(s).
$$
Now we recall several definitions and known facts related to the
$q$-Markov-WZ method (see \cite{moq}). An expression $H(n,k)$, in
the integer variables $n$ and $k,$ is called a {\it
$q$-hypergeometric term}  if the quotients
$$
\frac{H(n+1,k)}{H(n,k)} \qquad\mbox{and} \qquad
\frac{H(n,k+1)}{H(n,k)}
$$
are both rational functions with respect to $q^k$ and $q^n.$
 If for a given
$q$-hypergeometric term $H(n,k),$  there exists a polynomial
$P(n,k)$ in $q^k$ of the form
\begin{equation}
P(n,k)=a_0(n)+a_1(n)q^k+\cdots+a_L(n)q^{kL}, \label{p}
\end{equation}
for some non-negative integer $L,$ and P-recursive functions of
variable $q^n,$ $a_0(n), \ldots, a_L(n)$ such that
$F(n,k):=H(n,k)P(n,k)$ satisfies (\ref{WZ})
with some function $G,$ then a pair $(F,G)$ is called a {\it
$q$-Markov-WZ pair} associated with the kernel $H(n,k)$ ($q$-MWZ
pair for short). We call $G(n,k)$ a {\it $q$-MWZ mate} of
$F(n,k).$

\vspace{0.3cm}

{\bf Proposition.} \cite[Th.~7.2, 7.3]{moq} {\it For $q$-Markov-WZ
pair $(F,G),$ we have
$$
(i) \qquad
\sum_{k=0}^{\infty}F(0,k)-\lim_{L\to\infty}\sum_{k=0}^{\infty}F(L,k)
=\sum_{n=0}^{\infty}G(n,0)-\lim_{K\to\infty}\sum_{n=0}^{\infty}G(n,K),
$$
$$
(ii) \qquad
\sum_{k=0}^{\infty}F(0,k)=\sum_{n=0}^{\infty}(F(n,n)+G(n,n+1))-
\lim_{K\to\infty}\sum_{n=0}^{\infty}G(n,K),  \quad
$$
provided that both sides  converge.}

\vspace{0.2cm}

Now we can prove a $q$-analog of the Bailey-Borwein-Bradley
identity (\ref{eq04}).

\begin{theorem} \label{t1} Let $a$ be a complex number, with $|a|<1.$
Then for the generating function of the sequence
$\{\zeta[2k+2]\}_{k\ge 0},$ we have
$$
\sum_{k=0}^{\infty}\zeta[2k+2]a^{2k}=\sum_{n=1}^{\infty}
\frac{q^n}{[n]_q^2-a^2q^{2n}}= \sum_{n=1}^{\infty}
\frac{q^{n^2}(1+2q^n)}{\genfrac{[}{]}{0pt}{}{2n}{n}_q([n]_q^2-a^2q^{2n})}
\prod_{m=1}^{n-1}\frac{[m]_q^2-a^2(q^m+1)^2}{[m]_q^2-a^2q^{2m}}.
$$
In particular,
$$
\zeta[2]=\sum_{n=1}^{\infty}
\frac{q^{n^2}(1+2q^n)}{\genfrac{[}{]}{0pt}{}{2n}{n}_q[n]_q^2}.
$$
\end{theorem}

{\bf Proof.} We consider the $q$-hypergeometric term
\begin{equation}
H(n,k)=q^{k(2n+1)}\frac{(q\pm aq(1-q);q)_k}%
{(q\pm aq(1-q);q)_{n+k+1}} \label{h}
\end{equation}
and define the function $F(n,k)=H(n,k)P(n,k),$ where $P(n,k)$ is a
polynomial in $q^k$ of degree $L_1$ with unknown coefficients as
functions of $q^n.$ We are interested in finding a $q$-Markov-WZ
pair associated with $H(n,k).$ From the equality
$$
F(n+1,k)-F(n,k)=\frac{q^{k(2n+1)}(q\pm aq(1-q);q)_k}{(q\pm
aq(1-q);q)_{n+k+2}}P_1(n,k),
$$
where $P_1(n,k)$ is a polynomial in $q^k$ of degree $L_1+2,$ it
follows that we can determine a $q$-MWZ mate of $F(n,k)$ in the
form $G(n,k)=H(n,k)Q(n,k),$ where $Q(n,k)$ is a polynomial in
$q^k$ of degree $L_2$ with unknown coefficients as functions of
$q^n.$ Indeed, for such a choice we have
$$
G(n,k+1)-G(n,k)=\frac{q^{k(2n+1)}(q\pm aq(1-q);q)_k}{(q\pm
aq(1-q);q)_{n+k+2}}Q_1(n,k),
$$
where $Q_1(n,k)$ is a polynomial in $q^k$ such that $\deg
Q_1(n,k)=L_2+2,$ if $L_2\ne 1,$ and $\deg Q_1(n,k)=L_2+1,$ if
$L_2=1.$ Therefore, $(F,G)$ is a $q$-Markov-WZ pair iff
\begin{equation}
P_1(n,k)=Q_1(n,k) \quad \mbox{identically for all} \quad n,k.
\label{eq05}
\end{equation}
This implies that either $L_1=L_2$ and $L_2\ne 1,$ or $L_1=0,
L_2=1.$ On the other hand, equating coefficients of powers of
$q^k$ on both sides of (\ref{eq05}) we get a system of $L_1+3$
linear homogeneous equations with $L_1+L_2+2$ unknowns. In order
to guarantee a solution, we should at least have that
$L_1+L_2+2\ge L_1+3$ and hence $L_2\ge 1.$  We now show that there
is a non-zero solution of (\ref{eq05}) with the optimal choice
$L_2=1, L_1=0.$ To see this, define two functions
$$
F(n,k)=H(n,k)A(n), \quad G(n,k)=H(n,k)(B(n)+C(n)q^k),
$$
with $3$ unknown coefficients $A(n), B(n), C(n)$ as functions of
$q^n.$ Substituting $F, G$ into (\ref{WZ}) and dividing both sides
by $q^{k(2n+1)}\frac{(q\pm aq(1-q);q)_k}{(q\pm
aq(1-q);q)_{n+k+2}}$ we get that (\ref{WZ}) is equivalent to the
following equation of degree $2$ in the variable $q^k:$
\begin{equation}
\begin{split}
q^{2k}A(n+1)&-(1-q^{n+k+2}\pm aq^{n+k+2}(1-q))\cdot A(n)
=q^{2n+1}(1-q^{k+1}\pm aq^{k+1}(1-q)) \\ &\times (B(n)
+C(n)q^{k+1})-(1-q^{n+k+2}\pm aq^{n+k+2}(1-q))\cdot(B(n)+C(n)q^k).
\label{eq06}
\end{split}
\end{equation}
Equating all the coefficients of the powers of $q^k$ to zero in
(\ref{eq06}), we get a system of first order linear recurrence
equations with polynomial coefficients in $q^n$ for $A(n), B(n),
C(n):$
\begin{equation*}
\begin{split}
A(n)&=(1-q^{2n+1})B(n), \\
A(n)&=(1-q^n)B(n)-\frac{1-q^{2n+2}}{2q^{n+2}}C(n), \\
A(n+1)&=q^{2n+3}(1-a^2(1-q)^2) (qA(n)+(1-q)B(n))
+2q^{n+2}(1-q^{n+1}) C(n).
\end{split}
\end{equation*}
Resolving this system with respect to $A(n)$ we get
\begin{equation*}
\begin{split}
B(n)&=\frac{1}{1-q^{2n+1}}A(n), \\
C(n)&=\frac{-2q^{2n+2}}{(1+q^{n+1})(1-q^{2n+1})}A(n), \\
A(n+1)&=\frac{q^{2n+3}(1-q)^2[n+1]_q^2([n+1]_q^2-a^2(1+q^{n+1})^2)}%
{[2n+1]_q[2n+2]_q}A(n).
\end{split}
\end{equation*}
If we put $A(0)=q(1-q)^2,$ then
$$
A(n)=\frac{q^{(n+1)^2}(1-q)^{2n+2}}{\genfrac{[}{]}{0pt}{}{2n}{n}_q}
\prod_{m=1}^n([m]_q^2-a^2(1+q^m)^2),
$$
and a $q$-Markov-WZ pair associated with the kernel (\ref{h}) has
the form
\begin{equation*}
\begin{split}
F(n,k)&=\frac{q^{(n+1)^2+k(2n+1)}(1-q)^{2n+2}}{\genfrac{[}{]}{0pt}{}{2n}{n}_q}
\frac{(q\pm aq(1-q);q)_k}{(q\pm
aq(1-q);q)_{n+k+1}}\prod_{m=1}^n([m]_q^2-a^2(1+q^m)^2) \\[5pt]
G(n,k)&=\frac{q^{(n+1)^2+k(2n+1)}(1-q)^{2n+1}(1+q^{n+1}-2q^{k+2n+2})}%
{\genfrac{[}{]}{0pt}{}{2n+2}{n+1}_q[n+1]_q} \frac{(q\pm
aq(1-q);q)_k}{(q\pm aq(1-q);q)_{n+k+1}} \\
&\times\prod_{m=1}^n([m]_q^2-a^2(1+q^m)^2).
\end{split}
\end{equation*}
Now by the first summation formula from the Proposition we have
$$
\sum_{k=0}^{\infty}F(0,k)=\sum_{n=0}^{\infty}G(n,0)
$$
yielding the theorem. \qed

Applying the second summation formula to the pair $(F,G)$ we get a
$q$-analog of  \cite[Theorem 2]{he1}.

\begin{theorem} \label{t2}  Let $a$ be a complex number, with $|a|<1.$
Then for the generating function of the sequence
$\{\zeta[2k+2]\}_{k\ge 0},$ we have
$$
\sum_{k=0}^{\infty}\zeta[2k+2]a^{2k}=\sum_{n=1}^{\infty}
\frac{q^{3n^2-3n+1}p(n)}%
{\genfrac{[}{]}{0pt}{}{2n}{n}_q[n]_q([n]_q^2-a^2q^{2n})([2n]_q^2-a^2q^{4n})}
\prod_{m=1}^{n-1}\frac{[m]_q^2-a^2(q^m+1)^2}{[n+m]_q^2-a^2q^{2(n+m)}},
$$

\noindent where
\begin{equation*}
\begin{split}
p(n)&=(3q^{4n-1}+6q^{3n-1}-q^{2n}+8q^{2n-1}-4q^n+8q^{n-1}+1)[n]_q^3
-8q^n[n]_q^2 \\
&-a^2q^{4n-1}((3q^{2n}+4q^n-q+3)[n]_q-2).
\end{split}
\end{equation*}

\noindent In particular,
$$
\zeta[2]=\sum_{n=1}^{\infty}\frac{q^{3n^2-3n+1}((3q^{4n-1}+6q^{3n-1}-q^{2n}
+8q^{2n-1}-4q^n+8q^{n-1}+1)[n]_q-8q^n)}{\genfrac{[}{]}{0pt}{}{2n}{n}_q^3\,[n]_q^3}.
$$
\end{theorem}


\section{$q$-analogues of Markov's formula for $\zeta(3)$}
\label{S3}

In \cite{moq}, M.~Mohammed showed that for any $q$-proper
hypergeometric kernel $H(n,k),$ that is a function of the form
$$
H(n,k)=\frac{\prod_{i=1}^sQ(a_in+b_ik,c_i)}{\prod_{j=1}^tQ(u_jn+v_jk,w_j)}
q^{an^2+bnk+ck^2+dk+en}\xi^k,
$$
where $Q(m,c)=(1-cq)(1-cq^2)\cdots (1-cq^m),$ there exists a
non-negative integer $L$ and a polynomial $P(n,k)$ in the form
(\ref{p}) such that $F(n,k)=H(n,k)P(n,k)$ has a $q$-MWZ mate
$G(n,k)=F(n,k)Q(n,k),$ where $Q(n,k)$ is a ratio of two
$P$-recursive functions in $(q^n,q^k).$

Paper \cite{moq} is accompanied by the Maple package $q$-MarkovWZ
which, for a given $H(n,k),$ outputs the desired
$P(n,k)=\sum_{i=0}^La_i(n)q^{ki}$ and $G(n,k).$ Actually, the
output consists of a transition matrix $A$ between the $a_i$'s at
$n$ and the $a_i$'s at $n+1$
$$
(a_0(n+1), \ldots, a_L(n+1))^T=A(a_0(n), \ldots, a_L(n))^T,
$$
which denotes the $a_i(n)$'s subject to the initial conditions,
and a vector of  $b_i$'s, rational functions in the variables
$q^n$ and $q^k,$ whose dot product with the $a_i(n)$'s multiplied
by $H$ gives the $q$-MWZ mate
$G(n,k)=H(n,k)\sum_{i=0}^La_i(n)b_i(n,k).$

The Maple package $q$-MarkovWZ works only with a $q$-proper
hypergeometric kernel of the form
$$
H(n,k)=\frac{\prod_{i=1}^s(q;q)_{a_in+b_ik+c_i}}%
{\prod_{j=1}^t(q;q)_{u_jn+v_jk+w_j}}q^{\alpha_1n^2+
\alpha_2n+\alpha_3nk+\alpha_4}q^{\beta_1k^2+\beta_2k+\beta_3}\xi^k,
$$
where $a_i,b_i,c_i,u_i,v_i,w_i,\alpha_i,\beta_i$ are non-negative
specific integers and $\xi\in {\mathbb R}.$ This means that
identities of Theorems \ref{t1}, \ref{t2}, for example,  cannot be
obtained by running the Maple package $q$-MarkovWZ. However, it is
always possible to apply it to a pure $q$-series that does not
contain any additional parameters.

We consider a $q$-analog of $\zeta(3)$
$$
\zeta[3]=\sum_{n=1}^{\infty}\frac{q^{2n}}{[n]_q^3}, \qquad\qquad
|q|<1,
$$
or changing  $q$ by $q^{-1}$ for  simplicity in the sequel
$$
\zeta[3]=\frac{1}{q^3}\sum_{n=1}^{\infty}\frac{{ q}^n}{[n]_{
q}^3}, \qquad\qquad |{ q}|>1.
$$
We take the kernel
$$
H(n,k)=\frac{({ q};{ q})_k{ q}^k}%
{({ q};{ q})_{2n+k+1}(1-{ q}^{n+k+1})^2}, \qquad |q|>1,
$$
and corresponding to this kernel determine a $q$-Markov-WZ pair
$(F,G).$ Applying the Maple package $q$-MarkovWZ \cite{moq}, we
see that $F(n,k)$ is of the form
$F(n,k)=H(n,k)(a_0(n)+a_1(n)q^k),$ where
\begin{equation}
\left( \begin{array}{c} a_0(n+1) \\
a_1(n+1)
\end{array}
\right) = \left(
\begin{array}{cc}
0 & (q^{n+1}-1)^2q^{-2n-1} \\
-q^2(q^{n+1}-1)^2 & -2(q^{n+1}-1)^2q^{1-n}
\end{array}
\right) \left(
\begin{array}{c}
a_0(n) \\
a_1(n)
\end{array}
\right) \label{eq07}
\end{equation}
and the corresponding $q$-MWZ mate of $F(n,k)$ is
$G(n,k)=H(n,k)(a_0(n)b_0(n,k)+a_1(n)b_1(n,k))$ with
$$
b_0(n,k)=\frac{{ q}^{4n+2k+4}-{ q}^{2n+k+2}-2{ q}^{n+k+1}+{
q}^k+1}{{ q}^k({ q}^{n+1}-1)({ q}^{n+1}+1)({ q}^{2n+k+2}-1)},
$$
\begin{equation*}
\begin{split}
&b_1(n,k)=(1+2q^{5n+k+4}+q^{4n+k+4}-4q^{4n+k+3}+3q^{2n+k+2}
-2q^{n+k+1}+2q^{5n+2k+4} \\ &+q^{4n+2k+4}
-2q^{6n+2k+5}-2q^{5n+2k+5} +q^{4n+2k+3}
+q^{8n+3k+7}-2q^{n+1}+2q^{3n+2} \\
&-q^{4n+3}-q^{6n+3k+5})/(q^{2n+k+2}(q^{2n+1}-1)
(q^{2n+2}-1)(q^{2n+k+2}-1)).
\end{split}
\end{equation*}
Considering solutions of (\ref{eq07}) with different initial
conditions we get a family of $q$-analogues of Markov's formula
(\ref{z3}) for $\zeta(3).$

Solving system (\ref{eq07}) with the initial conditions
$a_0(0)=1,$ $a_1(0)=0$ we find
\begin{equation}
\begin{split}
a_0(n)&=(-1)^{n-1}(n-1)(q;q)_n^2q^{-n(n-1)/2}, \,\,\qquad n\ge 1, \\
a_1(n)&=(-1)^nn(q;q)_n^2q^{2-(n-1)(n-2)/2}, \qquad\qquad n\ge 0.
\end{split}
\label{eq08}
\end{equation}
(Formulas (\ref{eq08}) can be proved by induction on $n.$)  It is
not hard to check that $(F,G)$ is indeed a $q$-MWZ pair, i.e.,
\begin{equation}
F(n+1,k)-F(n,k)=G(n,k+1)-G(n,k). \label{12.5}
\end{equation}
 Now by the first
part of the Proposition, we get
\begin{equation}
\sum_{k=0}^{\infty}F(0,k)=\sum_{n=0}^{\infty}G(n,0)=
\sum_{n=0}^{\infty}H(n,0)(a_0(n)b_0(n,0)+a_1(n)b_1(n,0)).
\label{eq09}
\end{equation}
Taking into account that
$$
b_0(n,0)=\frac{q^{2n+2}+2q^{n+1}+2}{(q^{n+1}+1)^2},
$$
$$
b_1(n,0)=\frac{q^{4n+4}+2q^{3n+3}+q^{2n+3}+2q^{n+1}+q^2-2q}%
{q^{2n+2}(q^{n+1}+1)^2}+\frac{(q-1)^2}{q^{2n+2}(q^{2n+2}-1)(q^{n+1}+1)^2}
$$
and simplifying (\ref{eq09}) we get
\begin{equation*}
\begin{split}
& \sum_{k=0}^{\infty}\frac{q^k}{(1-q^{k+1})^3}=\sum_{n=0}^{\infty}
\frac{(-1)^n q^{-n(n-1)/2}}{\genfrac{[}{]}{0pt}{}{2n+2}{n+1}_q
(1-q^{n+1})^3(1+q^{n+1})}\left(q^{2n+2}+2q^{n+1}+2 \right.\\
&\left.+ n\Bigl(q^{n+1}(q^{n+1}-1)(q^{n+1}+2)+
\frac{q^{2n+3}-1}{q^{n+1}}+
\frac{(q-1)^2q^n}{q^{2n+1}-1}\Bigr)\right), \quad |q|>1.
\end{split}
\end{equation*}
Changing $q$ by $q^{-1}$ we have
\begin{equation}
\begin{split}
\zeta[3]&=\sum_{n=1}^{\infty}\frac{q^{2n}}{[n]_q^3}=
\sum_{n=1}^{\infty}\frac{(-1)^{n-1}q^{n(3n-1)/2}}%
{\genfrac{[}{]}{0pt}{}{2n}{n}_q[n]_q^3(1+q^n)}\left(q^n+2q^{2n}
+2q^{3n} \right.\\
&+\left.n(1-q)\Bigl((1+2q^n)[n]_q+q^{2n-1}[2n+1]_q+
\frac{q^{4n-2}}{[2n-1]_q}\Bigr)\right), \quad |q|<1.
\end{split}
\label{z31}
\end{equation}
Starting with the initial conditions $a_0(0)=0,$ $a_1(0)=1$ we get
the following solution of (\ref{eq07}):
$$
\begin{array}{cccc}
a_0(n) & = & (-1)^{n-1}n(q;q)_n^2q^{-1-n(n-1)/2}, \quad\qquad &
n\ge 0,
\\[3pt]
a_1(n) & = & (-1)^n(n+1)(q;q)_n^2q^{1-(n-1)(n-2)/2}, \quad\qquad &
n\ge 0.
\end{array}
$$
Then by (\ref{eq09}), after simplifying we get
\begin{equation*}
\begin{split}
&
\sum_{k=0}^{\infty}\frac{q^{2k}}{(1-q^{k+1})^3}=\sum_{n=0}^{\infty}
\frac{(-1)^n q^{1-n(n-1)/2}}{\genfrac{[}{]}{0pt}{}{2n+2}{n+1}_q
(1-q^{n+1})^3(1+q^{n+1})}\left(q^{3n+3}+2q^{2n+2}+q^{n+2}+2 \right.\\
&\left.+ nq^{n+1}(q^{n+1}-1)(q^{n+1}+2)+
\frac{nq^{2n+3}-n-1}{q^{n+1}}+
\frac{(n+1)(q-1)^2q^n}{q^{2n+1}-1}\right), \quad |q|>1,
\end{split}
\end{equation*}
or
\begin{equation}
\begin{split}
&\sum_{n=1}^{\infty}\frac{q^{n}}{[n]_q^3}=
\sum_{n=1}^{\infty}\frac{(-1)^{n-1}q^{n(3n-1)/2}}%
{\genfrac{[}{]}{0pt}{}{2n}{n}_q[n]_q^3(1+q^n)}\left(1+2q^n+q^{2n-1}
+2q^{3n}-q^{4n} \right.\\
&+\left.(1-q)\Bigl(n(1+2q^n)[n]_q+nq^{2n-1}[2n+1]_q+
\frac{(n+1)q^{4n-2}}{[2n-1]_q}\Bigr)\right), \quad |q|<1.
\end{split}
\label{z32}
\end{equation}
Note that any linear combination of the series figured on the left
of (\ref{z31}), (\ref{z32})
$$
\sum_{n=1}^{\infty}\frac{\alpha_1q^{2n}+\alpha_2q^n}{[n]_q^3},
\qquad |q|<1, \quad \alpha_1+\alpha_2\ne 0,
$$
can be considered as a $q$-analogue of $\zeta(3),$ since
$$
\lim_{q\to
1-}\sum_{n=1}^{\infty}\frac{\alpha_1q^{2n}+\alpha_2q^n}{[n]_q^3}
=(\alpha_1+\alpha_2)\zeta(3),
$$
and therefore summing formulas (\ref{z31}), (\ref{z32}) multiplied
by arbitrary constants we get an infinite family of different
$q$-analogues of Markov's formula (\ref{eq01}), among which
formula (\ref{z31}) has the simplest form.

\section{A $q$-analogue of Amdeberhan's series for $\zeta(3)$}

In this paragraph we obtain a $q$-analogue of the accelerated
series
$$
\zeta(3)=\frac{1}{4}\sum_{n=1}^{\infty}(-1)^{n-1}
\frac{56n^2-32n+5}{n^3(2n-1)^2\binom{2n}{n}\binom{3n}{n}}
$$
found by T.~Amdeberhan in \cite{am}. To do this, we apply the
second summation formula from the Proposition to the $q$-Markov-WZ
pair constructed in (\ref{12.5}). Then we get
\begin{equation*}
\begin{split}
&\sum_{k=1}^{\infty}\frac{q^{k-1}}{(1-q^{k})^3}=
\sum_{n=0}^{\infty}(F(n,n)+G(n,n+1))
=\sum_{n=0}^{\infty} \frac{(-1)^n(q;q)_n^3q^{-n-1-n(n-1)/2}}%
{(q;q)_{3n+3}(q^{2n+1}-1)^2(1+q^{n+1})^2} \\[3pt]
&\times (q^{2n+1}(nq^{2n+1}-n+1)(1-q^{3n+2})(1-q^{3n+3})(1+q^{n+1})^2 \\
&+
(q^{2n+1}-1)(nq^{8n+7}+nq^{7n+6}+q^{6n+5}-(2n-1)q^{5n+4}-(n+1)q^{4n+4}
\\ &-(n-2)q^{4n+3}-q^{3n+3}+(n+1)q^{3n+2}-2q^{2n+2}-q^{n+1}+n),
\qquad |q|>1,
\end{split}
\end{equation*}
or after changing $q$ by $q^{-1}$ we have
$$
\zeta[3]=\sum_{k=1}^{\infty}\frac{q^{2k}}{[k]_q^3}=
\sum_{n=1}^{\infty} \frac{(-1)^{n-1}q^{7n(n-1)/2+1}p(n)}%
{\genfrac{[}{]}{0pt}{}{2n}{n}_q\genfrac{[}{]}{0pt}{}{3n}{n}_q
[n]_q^3[2n-1]_q^2(1+q^n)^2}, \qquad |q|<1,
$$
where
\begin{equation*}
\begin{split}
p(n)&=[2n-1]_q^2q^{4n-1}(1+q^n+2q^{2n}+q^{3n})+[3n-1]_q[3n]_qq^{2n-1}(1+q^n)^2
\\
&+ (1-q)(n-1)(1+q^n)\left((1+q^n)[3n-1]_q[3n]_q \right.
\\ &-\left.q^{2n-1}[n]_q([n]_q+[2n]_q -[3n]_q-[4n-1]_q-[6n-1]_q)\right).
\end{split}
\end{equation*}

\section{Concluding remarks}

Note that while accelerated series for ordinary zeta values have a
geometric rate of convergence because of appearance of binomial
coefficients in denominators of the general terms, the accelerated
$q$-series for $q$-zeta values converge much faster in view of the
quantity $q^{an^2}$ containing in the $n$-th general form of the
series. The $q$-binomial coefficient
$\genfrac{[}{]}{0pt}{}{cn}{dn}_q$ doesn't influence  the speed of
convergence since from the asymptotic formula for $q$-factorial
(see \cite{man})
$$
[n]_q!=\frac{(q;q)_n}{(1-q)^n}=[2]_q^{1/2}\Gamma_{q^2}(1/2)
(1-q)^{-1/2-n}e^{\frac{\theta q^{n+1}}{(1-q)-q^{n+1}}}, \qquad
0<\theta <1,
$$
 it follows that $\genfrac{[}{]}{0pt}{}{cn}{dn}_q=O(1)$ as
$n\to\infty.$

Another remark is concerned with the fact that if $(F,G)$ is a
$q$-Markov-WZ pair, then $F_s(n,k):=F(sn,k)$ and
$G_s(n,k):=\sum_{i=0}^{s-1}G(sn+i,k),$ $s\in {\mathbb N},$ satisfy
(see \cite{az})
\begin{equation*}
F_s(n+1,k)-F_s(n,k)=G_s(n,k+1)-G_s(n,k)
\end{equation*}
and therefore we get the following summation formula:
\begin{equation}
\sum_{k=0}^{\infty}F(0,k)=\sum_{n=0}^{\infty}(F(sn,n)+\sum_{i=0}^{s-1}
G(sn+i,n+1))-\lim_{K\to\infty}\sum_{n=0}^{\infty}\sum_{i=0}^{s-1}
G(sn+i,K), \label{eqposl}
\end{equation}
which for $s=1$ coincides with the formula $(ii)$ from the
Proposition. By application of formula (\ref{eqposl}) with $s>1$
to the $q$-Markov-WZ pairs constructed in Sections \ref{S2},
\ref{S3}, one can derive additional  acceleration formulae for the
$q$-analogues of $\zeta(2)$ and $\zeta(3).$ These new $q$-series,
in general, converge faster as the parameter $s$ increases.


\begin{thebibliography}{99}
\bibitem{algr}
G. Almkvist, A. Granville, {\it Borwein and Bradley's Ap\'ery-like
formulae for $\zeta(4n+3),$} Experiment. Math., {\bf 8} (1999),
no. 2, 197-203.
\bibitem{am}
T. Amdeberhan,  {\it Faster and fasterconvergent  series for
$\zeta(3),$} Electron. J. Combinatorics {\bf 3(1)} (1996),
$\#$R13.
\bibitem{az}
T. Amdeberhan, D. Zeilberger, {\it Hypergeometric series
acceleration via the WZ method,} Elect. J. Combinatorics {\bf
4(2)} (1997), $\#$R3.
\bibitem{bb}
J. M.~Borwein, D. M.~Bradley, {\it Empirically determined
Ap\'ery-like formulae for $\zeta(4n+3)$,} Experiment. Math. {\bf
6} (1997), no. 3, 181-194.
\bibitem{bbb}
D. H. Bailey, J. M. Borwein, D. M. Bradley, {\it Experimental
determination of Ap\'ery-like identities for $zeta(2n+2),$}
Experiment. Math. {\bf 15} (2006), no. 3, 281-289.
\bibitem{br}
D. M. Bradley, {\it Hypergeometric functions related to series
acceleration formulas,} Contemporary Math. {\bf 457} (2008),
113-125.
\bibitem{he1}
Kh.~Hessami Pilehrood, T.~Hessami Pilehrood, {\it Generating
function identities for $\zeta(2n+2), \zeta(2n+3)$ via the
WZ-method,} Electron. J. Combinatorics {\bf 15} (2008), $\#$R35.
\bibitem{he2}
Kh.~Hessami Pilehrood, T.~Hessami Pilehrood, {\it Simultaneous
generation for zeta values by the Markov-WZ method,} DMTCS
(Discrete Math. Theor. Computer Science) {\bf 10:3} (2008),
115-124.
\bibitem{kkw}
M.~Kaneko, N.~Kurokawa, M.~Wakayama, {\it A variation of Euler's
approach to values of the Riemann zeta function,} Kyushu J. Math.
{\bf 57} (2003), 175-192.
\bibitem{ko}
M. Koecher, {\it Letter (German),} Math. Intelligencer, {\bf 2}
(1979/1980), no. 2, 62-64.
\bibitem{le}
D. Leshchiner, {\it Some new identities for $\zeta(k)$,} J. Number
Theory, {\bf 13} (1981), 355-362.
\bibitem{man}
M.~Mansour, {\it An asymptotic expansion of the $q$-gamma function
$\Gamma_q(x),$} J. Nonlinear Math. Physics {\bf 13} (2006), no. 4,
479-483.
\bibitem{ma}
A. A. Markoff, {\it M\'emoire sur la transformation des s\'eries
peu convergentes en s\'eries t\`res convergentes,} M\'emoires de
l'Academie Imp\'eriale des Sciences de St.-Petersbourg, VII
s\'erie, t. XXXVII, No.9 (1890).
\bibitem{moq}
M.~Mohammed, {\it The $q$-Markov-WZ method,} Annals of
Combinatorics {\bf 9} (2005), 205-221.
\bibitem{po} A.~van der Poorten, {\it A proof that Euler
missed... Ap\'ery's proof of the irrationality of $\zeta(3).$ An
informal report,} Math. Intelligencer {\bf 1} (1978/79), no. 4,
195-203.
\bibitem{ri}
 T. Rivoal, {\it Simultaneous generation of Koecher and
Almkvist-Grainville's Ap\'ery-like formulae,} Experiment. Math.,
{\bf 13} (2004), 503-508.
\end{thebibliography}
\end{document}